\documentclass[11pt]{amsart}
\usepackage[hmargin=2.5cm,vmargin=2.5cm]{geometry}
\usepackage{amsfonts}
\usepackage{amsthm}
\usepackage[english]{algorithm2e}
\usepackage{amsmath}
\usepackage{amssymb}
\usepackage{hyperref}
\usepackage[T1]{fontenc}
\newtheorem{Theorem}{Theorem}[section]
\newtheorem{Definition}[Theorem]{Definition}
\newtheorem{Remark}[Theorem]{Remark}
\newtheorem{Example}[Theorem]{Example}
\newtheorem{Lemma}[Theorem]{Lemma}
\newtheorem{Proposition}[Theorem]{Proposition}
\newtheorem{Corollary}[Theorem]{Corollary}

\title{Bimodules over Hom-Jordan and Hom-alternative algebras}
\author[Sylvain Attan]
       {Sylvain Attan
       \\\\
           D\'epartement de Math\'ematiques\\
       Universit\'e d'Abomey-calavi\\
       01 BP 4521, Cotonou 01, B\'enin\\
       syltane2010@yahoo.fr
       }
\begin{document}
\maketitle

\begin{abstract} In this paper,  bimodules over Hom-Jordan algebras and the ones over Hom-alternative algebras are defined. It is shown that bimodules over Jordan and alternative algebras 
are twisted into bimodules over Hom-Jordan and Hom-alternative algebras via endomorphisms respectively. Some relations between 
 bimodules over Hom-associative, Hom-Jordan and Hom-alternative algebras, are given.
\end{abstract}

\section{Introduction}
Algebras where the identities defining the structure are twisted by a homomorphism are called Hom-algebras. They have been intensively investigated in the literature recently.The theory of Hom-algebra started from Hom-Lie algebras introduced and discussed in
\cite{HAR1}, \cite {dlsds1}, \cite{dlsds2}, \cite {dlsds3}, motivated by quasi-deformations of Lie algebras of vector fields, in
particular q-deformations of Witt and Virasoro algebras. Hom-associative algebras were introduced in \cite{MAK3} while Hom-alternative and Hom-Jordan algebras are introduced in \cite{MAK1}\cite{YAU3} as twisted generalizations of alternative and Jordan algebra respectively.
The reader is referred to \cite{jtrmw} for applications of alternative algebras to projective geometry, buildings, and algebraic groups and to \cite{fgcht, pjjnew, so, tasfdv} for discussions about the important roles of Jordan algebras in physics, especially quantum mechanics.

The anti-commutator of a Hom-alternative algebra  gives rise to a Hom-Jordan algebra \cite{YAU3}. 
Star-ting with a Hom-alternative algebra $(A, \cdot,\alpha),$ it is known that the Jordan product
\begin{eqnarray}
x\ast y=\frac{1}{2}(x\cdot y+y\cdot x)\nonumber
\end{eqnarray}
gives a Hom-Jordan algebra $A^+=(A, \ast, \alpha).$ In other words, Hom-alternative algebras are Hom-Jordan-admissible \cite{YAU3}.

The notion of bimodule for a class of algebras defined by multilinear identities has been introduced by Eilenberg\cite{se}. If $\mathcal{H}$ is in the class of associative algebras or in the one of Lie algebras then this notion is the familiar one for which we are in possession of well-worked theories. The study of bimodule (or representation) of Jordan algebras was initiated by N. Jacobson \cite{Jacob1}.  Subsequently the alternative case was considered by Schafer \cite{Schaf1}.

Modules over an ordinary algebra has been extended to the ones of Hom-algebras in many works \cite{ibbm}, \cite{YS}, \cite{YAU4}, \cite{YAU2}. 

The aim of this paper is to introduce bimodules over Hom-alternative algebras and Hom-Jordan 
algebras and to discuss of some results about theses concepts. The paper is organized as follows. In section one, we recall basic notions related to Hom-algebras and modules over Hom-associative  algebras. Section two is devoted to the introduction of bimodules over Hom-alternative algebras.   Proposition \ref{sma} shows that from a given bimodule over a Hom-alternative algebra, a sequence of this kind of bimodules can be obtained. Theorem \ref{mamHa} shows that, a bimodule over alternative algebra gives rise to  a bimodule over the corresponding twisted algebra. It is also proved that a direct sum of a Hom-alternative algebra and a module over this Hom-algebra is again a Hom-alternative algebra (Theorem \ref{smHa} ).
In section three we introduce modules over Hom-Jordan algebras and prove similar results as in the previous section. Furthermore, it is proved that a left and right special module over a Hom-Jordan algebra with an additional  condition, has a bimodule structure over this Hom-algebra (Theorem \ref{ModPlus}). Finally Proposition \ref{bhamHJ} shows that a bimodule over a Hom-associative algebra has bimodule structure over its plus Hom-algebra. \\
All vector spaces are assumed to be  over a fixed ground field $\mathbb{K}$ of characteristic $0.$
\section{Preliminaries and main results}
We recall some basic notions, introduced in \cite{HAR1}, \cite{MAK3}, \cite{YAU4} related to Hom-algebras and while dealing of any binary operation  we will use juxtaposition in order to reduce the 
number of braces i.e., e.g., for $"\cdot ", $  $xy\cdot\alpha(z)$ means $(x\cdot y)\cdot\alpha(z).$ 
Also for the map $\mu: A^{\otimes 2}\longrightarrow A,$ we will sometimes $\mu(a\otimes b)$ as $\mu(a,b)$ or $ab$  for $a,b\in A$ and if $V$ is another vector 
space, $\tau_1: A\otimes V\longrightarrow V\otimes A$ (resp. $\tau_2: V\otimes A\longrightarrow A\otimes V$) denote the twist isomorphism $\tau_1(a\otimes v)=v\otimes a$ (resp.
$\tau_2(v\otimes a)=a\otimes v$).
\begin{Definition}
A Hom-module is a pair $(M,\alpha_M)$ consisting of a $\mathbb{K}$-module $M$ and 
a linear self-map $\alpha_M: M\longrightarrow M.$ A morphism 
$f: (M,\alpha_M)\longrightarrow (N,\alpha_N)$ of Hom-modules is a linear map 
 $f: M\longrightarrow N$ such that $f\circ\alpha_M=\alpha_N\circ f.$
\end{Definition}
\begin{Definition} (\cite{MAK3}, \cite{YAU4}) A Hom-algebra is a triple $(A,\mu,\alpha)$ in which $(A,\alpha)$ is a Hom-module, $\mu: A^{\otimes 2}\longrightarrow A$ i a linear map.
The Hom-algebra $(A,\mu,\alpha)$ is said to be  multiplicative if $\alpha\circ\mu=\mu\circ\alpha^{\otimes 2}$ (multiplicativity). A morphism 
$f: (A,\mu_A,\alpha_A)\longrightarrow (B,\mu_B,\alpha_B)$ of Hom-algebras is a morphism of the underlying Hom-modules such that and $f\circ\mu_A=\mu_B\circ f^{\otimes 2}.$
\end{Definition}
In this paper, we will only consider multiplicative Hom-algebras.
\begin{Definition}\label{DefRef}
Let $(A,\mu,\alpha)$ be a Hom-algebra.\\
(1) The Hom-associator of $A$ is the linear map $as_A: A^{\otimes 3}\longrightarrow A$  defined as $as_A=\mu\circ(\mu\otimes\alpha-\alpha\otimes\mu).$ A multiplicative Hom-algebra $(A,\mu,\alpha)$ is said to be Hom-associative algebra if $as_A=0$\\
(2) A  Hom-alternative algebra \cite{MAK1} is a multiplicative Hom-algebra $(A,\mu,\alpha)$ that satisfies both
\begin{eqnarray}
as_A(x,x,y)&=&0  \mbox{ (left Hom-alternativity)}\label{lhi} \mbox { and }\\
%\end{eqnarray}
 %\begin{eqnarray}
as_A(x,y,y)&=&0  \mbox{ (right Hom-alternativity)}\label{rhi}
\end{eqnarray}
for all $x,y\in A.$
\end{Definition}
\begin{Definition}\cite{YAU3}
A Hom-Jordan algebra is a multiplicative Hom-algebra $(A,\mu,\alpha)$ such that $\mu\circ\tau=\mu$
(commutativity of $\mu$) and the so-called Hom-Jordan identity 
\begin{eqnarray}
as_A(\mu(x,x,),\alpha(y),\alpha(x))=0 \label{idJord}
\end{eqnarray}
holds for all $x,y\in A.$
\end{Definition}
\begin{Remark}
In \cite{MAK1} Makhlouf defined a Hom-Jordan algebra as a commutative multiplicative
Hom-algebra satisfying $as_A (x^2,y,\alpha(x))=0,$ which becomes the identity (\ref{idJord}) if $y$ is replaced by $\alpha(y).$
\end{Remark}
%It is proved in\cite{•} that the Hom-Jordan identity (\ref{idJord}) is equivalent to the linearised %Hom-Jordan identity
%\begin{eqnarray}
%\circlearrowright_{x,y,z}as_A(\alpha(x),\alpha(w),\mu(y,z))=0
%\end{eqnarray}
The proof of the following result can be found in \cite{YAU3} where the product $\ast$  there, differs from the one given here  by a factor of $\frac{1}{2}$.
\begin{Proposition}
Let $(A,\mu,\alpha)$ be a Hom-alternative algebra. Then $A^+=(A,\ast,\alpha)$ is a Hom-Jordan algebra where 
$x\ast y=xy+yx$ for all $x,y\in A.$
\end{Proposition}
 A. Makhlouf proved that the plus algebra of any Hom-associative algebra is a Hom-Jordan 
algebra as defined in \cite{MAK1}. Here, we prove the same result for the Hom-Jordan algebra as defined in \cite{YAU3} (see also Definition \ref{DefRef}  above). 
\begin{Proposition}
Let $(A,\cdot, \alpha)$ be a Hom-associative algebra. Then $A^+=(A,\ast,\alpha)$ is a Hom-Jordan algebra where 
$x\ast y=xy+yx$ for all $x,y\in A.$
\end{Proposition}
{\bf Proof.} The commutativity of $\ast$ is obvious. We compute the Hom-Jordan identity as follows:
\begin{eqnarray}
as_{A^+}(x^2,\alpha(x),\alpha(y))&=&(x^2\ast\alpha(y))\ast\alpha^2(x)-\alpha(x^2)\ast(\alpha(y)\ast\alpha(x))\nonumber\\
&=&(x^2\cdot\alpha(y))\cdot\alpha^2(x)+(\alpha(y)\cdot x^2)\cdot\alpha^2(x)+
\alpha^2(x)\cdot(x^2\cdot\alpha(y)\nonumber\\
&&+\alpha^2(x)\cdot(\alpha(y)\cdot x^2)-\alpha(x^2)\cdot(\alpha(y)\cdot\alpha(x))-\alpha(x^2)\cdot(\alpha(x)\cdot\alpha(y))\nonumber\\
&&-(\alpha(y)\cdot\alpha(x))\cdot\alpha(x^2)-(\alpha(x)\cdot\alpha(y))\cdot\alpha(x^2) \mbox{(by a direct computation)}\nonumber\\
&=&(\alpha(y)\cdot x^2)\cdot\alpha^2(x)+
\alpha^2(x)\cdot(x^2\cdot\alpha(y)-\alpha(x^2)\cdot(\alpha(x)\cdot\alpha(y))\nonumber\\
&&-(\alpha(y)\cdot\alpha(x))\cdot\alpha(x^2) \mbox{(by the Hom-associativity)}\nonumber\\
&=&(\alpha(y)\cdot x^2)\cdot\alpha^2(x)+
\alpha^2(x)\cdot(x^2\cdot\alpha(y)-(\alpha(x)\cdot\alpha(x))\cdot\alpha(x\cdot y)\nonumber\\
&&-\alpha(yx)\cdot(\alpha(x)\cdot\alpha(x)) \mbox{(by the multiplicativity)}\nonumber\\
&=&(\alpha(y)\cdot x^2)\cdot\alpha^2(x)+
\alpha^2(x)\cdot(x^2\cdot\alpha(y)-\alpha^2(x)\cdot(\alpha(x)\cdot(x\cdot y))\nonumber\\
&&-((yx)\cdot\alpha(x))\cdot\alpha^2(x)\mbox{(by the the Hom-associativity)}\nonumber\\
&=&0 \mbox{(by the the Hom-associativity)}\nonumber
\end{eqnarray}
Then $A^+=(A,\ast,\alpha)$ is a Hom-Jordan algebra.\hfill $\square$\\
\\
Let us finish this section by the following definitions which will be used in next sections.
\begin{Definition} Let $(A,\mu,\alpha_A)$ be any Hom-algebra and $(V,\alpha_V)$ be a Hom-module.
\begin{enumerate}
\item A left (resp. right) structure map on $V$ is a morphism $\rho_l: A\otimes V\longrightarrow V,$
$a\otimes v\longmapsto a\cdot v$ (resp. $\rho_r: V\otimes A\longrightarrow V,$
$v\otimes a\longmapsto v\cdot a$) of Hom-modules.
\item Let $\rho_l$ and $\rho_r$ be structure maps on $V.$ Then the module Hom-associator of $V$ is a trilinear map $as_{A,V}$ defined as:
\begin{eqnarray}
as_{A,V}\circ Id_{V\otimes A\otimes A}=\rho_r\circ(\rho_r\otimes\alpha_A)-
\rho_l\circ(\alpha_V\otimes\mu)\nonumber\\
as_{A,V}\circ Id_{A\otimes V\otimes A}=\rho_r\circ(\rho_l\otimes\alpha_A)-
\rho_l\circ(\alpha_A\otimes\rho_r)\nonumber\\
as_{A,V}\circ Id_{A\otimes A\otimes V}=\rho_l\circ(\mu\otimes\alpha_V)-
\rho_l\circ(\alpha_A\otimes\rho_l)\nonumber
\end{eqnarray}
%that is
%\begin{eqnarray}
%as_V(a,v,b)=(a\cdot v)\cdot\alpha_A(b)-\alpha_A(a)\cdot(v\cdot b)\label{mix1}\nonumber\\
%as_V(v,a,b)=(v\cdot a)\cdot\alpha_A(b)-\alpha_V(v)\cdot(ab)\label{mix2}\nonumber\\
%as_V(a,b,v)=(ab)\cdot\alpha_V(v)-\alpha_A(a)\cdot(b\cdot v)\label{mix3}\nonumber
%\end{eqnarray}
%for all $a,b\in A$ and $v\in V$
\end{enumerate}
\end{Definition}
\begin{Remark}
The module Hom-associator given above is a generalization of the one given in \cite{ibbm}.
\end{Remark}
\begin{Definition} \cite{YAU4, YAU2} 
Let $(A,\mu,\alpha_A)$ be a Hom-associative algebra and $(M, \alpha_M)$ be a Hom-module.
\begin{enumerate}
\item A (left) $A$-module structure on $M$ consists of a morphism $\rho: A\otimes M\longrightarrow M$ of Hom-modules, such that
\begin{eqnarray}
\rho\circ(\alpha_A\otimes\rho)=\rho\circ(\mu\otimes\alpha_M)\label{LeftAssMod}
\end{eqnarray}
\item A right $A$-module structure on $M$ consists of a morphism $\rho: M\otimes A\longrightarrow M$ of Hom-modules, such that
\begin{eqnarray}
\rho\circ(\alpha_M\otimes\mu)=\rho\circ(\rho\otimes\alpha_A)\label{RightAssMod}
\end{eqnarray}
\item An $A$-bimodule structure on $M$ consists of two structure maps $\rho_l: A\otimes M\longrightarrow M$ and $\rho_r: M\otimes A\longrightarrow M$ such that $(M, \alpha_M,\rho_l)$ is 
a left $A$-module, $(M, \alpha_M,\rho_r)$ is a right $A$-module and that the following Hom-associativity (or operator commutativity) condition holds: 
\begin{eqnarray}
\rho_l\circ(\alpha_A\otimes\rho_r)=\rho_r\circ(\rho_l\otimes\alpha_A) \label{Assoc1}
\end{eqnarray}
\end{enumerate}
\end{Definition} 
\section{Hom-bimodule over Hom-alternative algebras}
In this section, we give the definition of modules over  Hom-alternative algebras. We prove that from a given bimodule over a Hom-alternative algebra, a sequence of this kind of bimodules can be constructed. It is also proved that a direct sum of a Hom-alternative algebra 
and a bimodule over this Hom-algebra, is a Hom-alternative algebra called a split null extension of the considered Hom-algebra.

 First, we start by the following notion, due to \cite{ibbm}, where it is called a module over a left (resp. right) Hom-alternative algebra. However, we call it a  left (resp. right) module over a Hom-alternative algebra in this paper to unify our terminologies.
\begin{Definition}
Let $(A,\mu,\alpha_A)$ be a Hom-alternative algebra.
\begin{enumerate}
\item A left $A$-module is a Hom-module $(V,\alpha_V)$ that comes equipped with a structure map $\rho_l:A\otimes V\longrightarrow V$ ($\rho_l(a\otimes v)=a\cdot v$)  such that
\begin{eqnarray}
as_{A,V}(x,y,v)=-as_{A,V}(y,x,v)
\end{eqnarray}
for all $x, y\in A$ and $v\in V.$
\item A right $A$-module is a Hom-module $(V,\alpha_V)$ that comes equipped with a structure map $\rho_r:V\otimes A\longrightarrow V$ ($\rho_r(v\otimes a)=v\cdot a$)   such that
\begin{eqnarray}
as_{A,V}(v,x,y)=-as_{A,V}(v,y,x)
\end{eqnarray}
for all $x, y\in A$ and $v\in V.$
\end{enumerate}
\end{Definition}
Now, as a generalization of  bimodules over alternative algebras, one has:
\begin{Definition}
Let $(A,\mu,\alpha_A)$ be a Hom-alternative algebra.\\
(i) An $A$-bimodule is a Hom-module $(V,\alpha_V)$ that comes equipped with a (left) structure map  $\rho_l:A\otimes V\longrightarrow V$ ($\rho_l(a\otimes v)=a\cdot v$) and a (right) structure map $\rho_r:V\otimes A\longrightarrow V$ ($\rho_r(v\otimes a)=v\cdot a$) 
such that the following equalities:
%\begin{eqnarray}
%\alpha_V\circ\rho_l=\rho_l\circ(\alpha_A\otimes\alpha_V)\label{ra1}\\
%\alpha_V\circ\rho_r=\rho_l\circ(\alpha_V\otimes\alpha_A)\label{ra2}
%\end{eqnarray}
\begin{eqnarray}
as_{A,V}(a,v,b)=-as_{A,V}(v,a,b)=as_{A,V}(b,a,v)=-as_{A,V}(a,b,v)\label{ra3}
\end{eqnarray}
hold for all $a,b,c\in A$ and $v\in V.$\\
(ii) A morphism $f :(V, \alpha_V )\longrightarrow (W, \alpha_W )$ of $A$-bimodules is a morphism of the underlying Hom-modules such that
\begin{eqnarray}
f\circ\rho_l=\rho_l\circ(Id_A\otimes f)\nonumber\\
f\circ\rho_r=\rho_l\circ(f \otimes Id_A)\nonumber
\end{eqnarray}
\end{Definition}
\begin{Remark}\label{rem2}
Since the field is of characteristic $0$, the relation (\ref{ra3}) is 
equivalent to\\
$\left\{
\begin{array}{rl}
as_{A,V}(a,v,b)=-as_{A,V}(v,a,b)=as_V(b,a,v)\\
as_{A,V}(a,a,v)=0
\end{array}
\right.$
\end{Remark}
\begin{Example} Here are some examples of $A$-modules.\\
(1)Let $(A,\mu,\alpha_A)$ be a Hom-alternative algebra. Then $(A,\alpha_A)$ is an $A$-bimodule where the structure maps are $\rho_l(a,b)=\mu(a,b)$ and $\rho_r(a,b)=\mu(a,b).$\\
(2) If $(A,\mu)$ is an alternative algebra and $M$ is an $A$-bimodule \cite{Jacob2} in the usual sense then $(M,Id_M)$ is an $\mathbb{A}$-bimodule where $\mathbb{A}=(A,\mu, Id_A)$ is a Hom-alternative algebra.
\end{Example}
The following result describes a sequence of bimodules over a Hom-alternative algebra by twisting 
the structure maps of a given module over this Hom-algebra. 
\begin{Proposition}\label{sma}
Let $(A,\mu,\alpha_A)$ be a Hom-alternative algebra and  $(V,\alpha_V)$  be an $A$-bimodule with the structure maps $\rho_l$ and $\rho_r$. Then the maps
\begin{eqnarray}
\rho_l^{(n)}=\rho_l\circ(\alpha_A^n\otimes Id_V)\label{nma1}\\
\rho_r^{(n)}=\rho_r\circ(Id_V\otimes \alpha_A^n)\label{nma2}
\end{eqnarray} 
give the Hom-module $(V,\alpha_V)$ the structure of an $A$-bimodule that we denote 
by $V^{(n)}$
\end{Proposition}
{\bf Proof.}
It is clear that $\rho_l^{(n)}$  and $\rho_r^{(n)}$ are structure maps on $V^{(n)}.$ Next, if we observe that for all $x,y\in A$ and $v\in V,$
\begin{eqnarray}
as_{A,V^{(n)}}(x,v,y)&=&\rho_r^n(\rho_l^n(x,v),\alpha_A(y))-\rho_l^n(\alpha_A(x),
\rho_r^n(v,y))\nonumber\\
&=&\rho_r(\rho_l(\alpha_A^n(x),v),\alpha_A^{n+1}(y))-\rho_l(\alpha_A^{n+1}(x),
\rho_r(v,\alpha_A^n(y))\nonumber\\
&=&as_{A,V}(\alpha_A^n(x),v,\alpha_A^n(y))\nonumber
\end{eqnarray}
and similarly
\begin{eqnarray}
as_{A,V^{(n)}}(v,x,y)&=&as_{A,V}(v,\alpha_A^n(x),\alpha_A^n(y))\nonumber\\
as_{A,V^{(n)}}(y,x,v)&=&as_{A,V}(\alpha_A^n(y),\alpha_A^n(x),v)\nonumber\\
as_{A,V^{(n)}}(x,y,v)&=&as_{A,V}(\alpha_A^n(x),\alpha_A^n(y),v)\nonumber
\end{eqnarray}
then (\ref{ra3}) for $V^{(n)}$ follows from the one in $V.$\hfill $\square$\\

 We know that alternative algebras can be deformed into Hom-alternative algebras via an endomorphism. The following result shows that bimodules over alternative algebras can  be deformed into bimodules over Hom-alternative algebras via an endomorphism. This result provides
a large class of examples of modules over Hom-alternative algebras
\begin{Theorem}\label{mamHa}
Let $(A,\mu)$ be an alternative algebra, $V$ be a $A$-bimodule with the structure maps
$\rho_l$ and $\rho_r$, $\alpha_A$ be an endomorphism of the alternative algebra $A$ 
and $\alpha_V$ be a linear self-map of $V$ such that $\alpha_V\circ\rho_l=\rho_l\circ(\alpha_A\otimes\alpha_V)$ and 
$\alpha_V\circ\rho_r=\rho_r\circ(\alpha_V\otimes\alpha_A)$\\
Write $A_{\alpha_A}$ for the Hom-alternative algebra $(A,\mu_{\alpha_A},\alpha_A)$ and 
$V_{\alpha_V}$ for the Hom-module $(V,\alpha_V).$ Then the maps:
\begin{eqnarray}
\tilde{\rho_l}=\alpha_V\circ\rho_l \mbox{ and }
\tilde{\rho_r}=\alpha_V\circ\rho_r
\end{eqnarray}
give the Hom-module $V_{\alpha_V}$ the structure of an $A_{\alpha_A}$-bimodule.
\end{Theorem}
{\bf Proof.} Trivially, $\tilde{\rho_l}$  and $\tilde{\rho_l}$ are structure maps on $V_{\alpha_V}.$ The proof of (\ref{ra3}) for $V_{\alpha_V}$ follows directly by  the fact that $as_{A,V_{\alpha_V}}=\alpha_V^2\circ as_{A,V}$ and the relation (\ref{ra3}) in $V.$\hfill $\square$
\begin{Corollary}
Let $(A,\mu)$ be an alternative algebra, $V$ be a $A$-bimodule with the structure maps
$\rho_l$ and $\rho_r$, $\alpha_A$ an endomorphism of the alternative algebra $A$ 
and $\alpha_V$ be a linear self-map of $V$ such that $\alpha_V\circ\rho_l=\rho_l\circ(\alpha_A\otimes\alpha_V)$ and 
$\alpha_V\circ\rho_r=\rho_r\circ(\alpha_V\otimes\alpha_A)$\\
Write $A_{\alpha_A}$ for the Hom-alternative algebra $(A,\mu_{\alpha_A},\alpha_A)$ and 
$V_{\alpha_V}$ for the Hom-module $(V,\alpha_V).$ Then the maps:
\begin{eqnarray}
\tilde{\rho_l^n}=\rho_l\circ(\alpha_A^{n+1}\otimes\alpha_V) \mbox{ and }
\tilde{\rho_r^n}=\rho_l\circ(\alpha_V\otimes\alpha_A^{n+1})
\end{eqnarray}
give the Hom-module $V_{\alpha}$ the structure of an $A_{\alpha_A}$-bimodule for each $n\in\mathbb{N}$. \hfill $\square$
\end{Corollary}
The following lemma will be used to prove the next theorem.
\begin{Lemma}
Let $(A,\mu,\alpha_A)$ be a Hom-alternative algebra and  $(V,\alpha_V)$  be an $A$-bimodule with the structure maps $\rho_l$ and $\rho_r$. Then the following relation
\begin{eqnarray}
as_{A,V}(v, a,a)=0 \label{lem2}
\end{eqnarray}
holds for all $a\in A$ and $v\in V.$
\end{Lemma}
{\bf Proof.} From the second and the first expression and the second and the third
expression of the relation (\ref{ra3}), we have  respectively for all $a,b\in V$ and $v\in V,$ $-as_{A,V}(v,a,b)=as_{A,V}(a,v,b)$ and $as_{A,V}(v,b,a)=-as_{A,V}(a,b,v).$ If we 
observe that from the first and the last expression of (\ref{ra3}),
$as_{A,V}(a,v,b)=-as_{A,V}(a,b,v),$ we then get $-as_{A,V}(v,a,b)=as_{A,V}(v,b,a).$ The later equality is equivalent to $as_{A,V}(v,a,a)=0$ since the field $\mathbb{K}$ is of characteristic  $0.$ \hfill $\square$\\

The following result shows that a direct sum of a Hom-alternative algebra and a bimodule over this Hom-algebra, is still a Hom-alternative, called the split null extension determined by the
given bimodule.
\begin{Theorem}\label{smHa}
Let $(A,\mu,\alpha_A)$ be a Hom-alternative algebra and  $(V,\alpha_V)$  be an $A$-bimodule with the structure maps $\rho_l$ and $\rho_r$. If define on $A\oplus V$ the bilinear maps
$\tilde{\mu}: (A\oplus V)^{\otimes 2}\longrightarrow A\oplus V,$ 
$\tilde{\mu}(a+m,b+n):=ab+a\cdot n+m\cdot b$  and the linear map
$\tilde{\alpha}: A\oplus V\longrightarrow A\oplus V,$ $\tilde{\alpha}(a+m):=\alpha_A(a)+\alpha_V(m),$ then $E=(A\oplus V,\tilde{\mu},\tilde{\alpha})$ is 
a Hom-alternative algebra.
\end{Theorem}
{\bf Proof.} The multiplicativity of $\tilde{\alpha}$ with respect to $\tilde{\mu}$ follows from the one of $\alpha$ with respect to $\mu$ and the that that $\rho_l$ and $\rho_r$ are morphisms of Hom-modules. Next
\begin{eqnarray}
as_E(a+m,a+m,b+n)&=&\tilde{\mu}(\tilde{\mu}(a+m,a+m),\tilde{\alpha}(b+n))
-\tilde{\mu}(\tilde{\alpha}(a+m),\tilde{\mu}(a+m,b+n))\nonumber\\
&=&\tilde{\mu}(a^2+a\cdot m+m\cdot a,\alpha_A(b)+\alpha_V(n))\nonumber\\
&&-\tilde{\mu}(\alpha_A(a)+\alpha_V(m),ab+a\cdot n+m\cdot b)\nonumber\\
&=&a^2\alpha_A(b)+a^2\cdot\alpha_V(n)+(a\cdot m)\cdot\alpha_A(b)+(m\cdot a)\cdot\alpha_A(b)\nonumber\\
&&-\alpha_A(a)(ab)-\alpha_A(a)\cdot(a\cdot n)-\alpha_A(a)\cdot(m\cdot b)-\alpha_V(m)\cdot(ab)\nonumber\\
&=&\underbrace{as_A(a,a,b)}_{0}+\underbrace{as_V(a,a,n)}_{0}+
\underbrace{as_{A,V}(a,m,b)+as_{A,V}(m,a,b)}_{0}\nonumber\\
  &&\mbox{ ( by (\ref{lhi}), Remark \ref{rem2} and (\ref{ra3}) )}\nonumber\\
&&=0\nonumber\\
%\end{eqnarray}
%\begin{eqnarray}
as_E(a+m,b+n,b+n)&=&\tilde{\mu}(\tilde{\mu}(a+m,b+n),\tilde{\alpha}(b+n))
-\tilde{\mu}(\tilde{\alpha}(a+m),\tilde{\mu}(b+n,b+n))\nonumber\\
&=&\tilde{\mu}(ab+a\cdot n+m\cdot b,\alpha_A(b)+\alpha_V(n))\nonumber\\
&&-\tilde{\mu}(\alpha_A(a)+\alpha_V(m),b^2+b\cdot n+b\cdot b)\nonumber\\
&=&(ab)\alpha_A(b)+(ab)\cdot\alpha_V(m)+(a\cdot n)\cdot\alpha_A(b)+
(m\cdot b)\cdot\alpha_A(b)\nonumber\\
&&-\alpha_A(a)(b^2)-\alpha_A(a)\cdot(b\cdot n)-\alpha_A(a)\cdot(n\cdot b)
-\alpha_V(m)\cdot b^2\nonumber\\
&=&\underbrace{as_A(a,b,b)}_{0}+\underbrace{as_{A,V}(a,b,n)+as_{A,V}(a,n,b)}_{0}
+\underbrace{as_{A,V}(m,b,b)}_{0}\nonumber\\
&&\mbox{ ( by (\ref{rhi}), (\ref{ra3}) and (\ref{lem2}) )}\nonumber\\
&=& 0\nonumber
\end{eqnarray}
We then conclude  that $(A\oplus V,\tilde{\mu},\tilde{\alpha})$ is a Hom-alternative algebra.
\hfill $\square$
\section{Hom-bimodule over Hom-Jordan algebras}
In this section, we define and study bimodules over Hom-Jordan algebras.  It is observed that 
similar results  for bimodules over Hom-alternative algebras, hold for the ones over Hom-Jordan algebras. Some of them, require an additional condition. Furthermore, on one hand, relations 
between bimodules over Hom-associative algebras and the ones over Hom-Jordan algebras are given and the other hand, relations  between left (resp. right) modules over Hom-alternative algebras and 
left(resp. right) special modules over Hom-Jordan algebras are proved. First  we have:
\begin{Definition}
Let $(A,\mu,\alpha_A)$ be a Hom-Jordan algebra. 
\begin{enumerate}
\item A right $A$-module is a Hom-module $(V,\alpha_V)$ that comes equipped with a right structure map $\rho_r:V\otimes A\longrightarrow V$ ($\rho_r(v\otimes a)=v\cdot a$) such that the following conditions:
%\begin{eqnarray}
%\alpha_V\circ\rho_r=\rho_l\circ(\alpha_V\otimes\alpha_A)\label{rc1}
%\end{eqnarray}
\begin{eqnarray}
&&\alpha_V(v\cdot a)\cdot\alpha_A(bc)+\alpha_V(v\cdot b)
\cdot\alpha_A(ca)+\alpha_V(v\cdot c)\cdot\alpha_A(ab)\nonumber
\end{eqnarray}
\begin{eqnarray}
&&=(\alpha_V(v)\cdot bc)\cdot\alpha_A^2(a)+(\alpha_V(v)\cdot ca)\cdot\alpha_A^2(b)+
(\alpha_V(v)\cdot ab)\cdot\alpha_A^2(c)\label{rc2}\\
&&\alpha_V(v\cdot a)\cdot\alpha_A(bc)+\alpha_V(v\cdot b)
\cdot\alpha_A(ca)+\alpha_V(v\cdot c)\cdot\alpha_A(ab)\nonumber\\
&&=((v\cdot a)\cdot\alpha_A(b))\cdot\alpha_A^2(c)+((v\cdot c)\cdot
\alpha_A(b))\cdot\alpha_A^2(a)+\alpha_V^2(v)\cdot((ac)\alpha_A(b))\label{rc3}
\end{eqnarray}
hold for all $a,b,c\in A$ and $v\in V.$\\
%(ii) A morphism $f :(V, \alpha_V )\longrightarrow (W, \alpha_W )$ of right Hom-$A$-modules is a %morphism of the underlying Hom-modules such that
%\begin{eqnarray}
%f\circ\rho_{r_V}=\rho_{r_W}\circ(f \otimes Id_A)\nonumber
%\end{eqnarray}
%\end{Definition}
%\begin{Remark}
%Let $(A,\mu,\alpha_A)$ be a Hom-Jordan algebra. \\
\item A left $A$-module is a Hom-module $(V,\alpha_V)$ that comes equipped with a left structure map $\rho_l:A\otimes V\longrightarrow V$ ($\rho_l(a\otimes v)=a\cdot v$) 
such that the following conditions:
%\begin{eqnarray}
% \alpha_V\circ\rho_l=\rho_l\circ(\alpha_A\otimes\alpha_V) \label{lc1}
%\end{eqnarray}
\begin{eqnarray}
&&\alpha_A(bc)\cdot\alpha_V(a\cdot v)+\alpha_A(ca)\cdot\alpha_V(b\cdot v)+\alpha_A(ab)\cdot\alpha_V(c\cdot v)\nonumber\\
&&=\alpha_A^2(a)\cdot(bc\cdot\alpha_V(v))+\alpha_A^2(b)\cdot(ca\cdot\alpha_V(v))+\alpha_A^2(c)\cdot(ab\cdot\alpha_V(v))
\label{lc2}\\
&&\alpha_A(bc)\cdot\alpha_V(a\cdot v)+\alpha_A(ca)\cdot\alpha_V(b\cdot v)+\alpha_A(ab)\cdot\alpha_V(c\cdot v)\nonumber\\
&&=\alpha_A^2(c)\cdot(\alpha_A(b)\cdot(a\cdot v))+\alpha_A^2(a)\cdot(\alpha_A(b)\cdot(c\cdot v))+((ac)\alpha_A(b))\cdot\alpha_V^2(v)\label{lc3}
\end{eqnarray}
hold for all $a,b,c\in A$ and $v\in V.$\\
%(ii) A morphism $f :(V, \alpha_V )\longrightarrow (W, \alpha_W )$ of left Hom-$A$-modules is a %morphism of the underlying Hom-modules such that
%\begin{eqnarray}
%f\circ\rho_{l_V}=\rho_{l_W}\circ(Id_A\otimes f)\nonumber
%\end{eqnarray}
%\end{Remark}
\end{enumerate}
\end{Definition}
The following result allows to introduce the notion of right special modules over Hom-Jordan algebras.
\begin{Theorem}\label{thmright}
Let $(A,\mu,\alpha_A)$ be a Hom-Jordan algebra, $(V,\alpha_V)$ be a Hom-module and $\rho_r: V\otimes A\rightarrow V,$ $(a\otimes v\mapsto v\cdot a),$ be a bilinear map satisfying
\begin{eqnarray}
\alpha_V\circ\rho_r=\rho_l\circ(\alpha_V\otimes\alpha_A) \label{rc5}\\
\alpha_V(v)\cdot(ab)=(v\cdot a)\cdot\alpha_A(b)+(v\cdot b)\cdot\alpha_A(a)\label{rc4}
\end{eqnarray}
for all $(a,b)\in A^2$ and $v\in V.$ Then $(V,\alpha, \rho_r)$ is a right $A$-module called
a right special $A$-module.
\end{Theorem}
{\bf Proof.} It suffices to prove $(\ref{rc2})$ and $(\ref{rc3}).$ For all $(a,b)\in A^2$ and $v\in V,$ we have:
\begin{eqnarray}
&&\alpha_V(v\cdot a)\cdot\alpha_A(bc)+\alpha_V(v\cdot b)
\cdot\alpha_A(ca)+\alpha_V(v\cdot c)\cdot\alpha_A(ab)\nonumber\\
&=&\alpha_V(v\cdot a)\cdot\alpha_A(b)\alpha_A(c)+\alpha_V(v\cdot b)
\cdot\alpha_A(c)\alpha_A(a)+\alpha_V(v\cdot c)\cdot\alpha_A(a)\alpha_A(b) (\mbox{ multiplicativity} )\nonumber\\
&=&((v\cdot a)\cdot\alpha_A(b))\cdot\alpha^2(c)+((v\cdot a)\cdot\alpha_A(c))\cdot\alpha^2(b)+
((v\cdot b)\cdot\alpha_A(c))\cdot\alpha^2(a)\nonumber\\
&&+((v\cdot b)\cdot\alpha_A(a))\cdot\alpha^2(c)+((v\cdot c)\cdot\alpha_A(a))\cdot\alpha^2(b)+((v\cdot c)\cdot\alpha_A(b))\cdot\alpha^2(a) \mbox{ (by (\ref{rc4}) )} \nonumber\\
&=&[\alpha_V(v)\cdot ab-(v\cdot b)\alpha_A(a)]\cdot\alpha^2(c)+((v\cdot a)\cdot\alpha_A(c))\cdot\alpha^2(b)\nonumber\\
&&+[\alpha_V(v)\cdot bc-(v\cdot c)\alpha_A(b)]\cdot\alpha^2(a)+((v\cdot b)\cdot\alpha_A(a))\cdot\alpha^2(c)\nonumber\\
&&+[\alpha_V(v)\cdot ca-(v\cdot a)\cdot\alpha_A(c)]\cdot\alpha^2(b)+((v\cdot c)\cdot\alpha_A(b))\cdot\alpha^2(a) \mbox{ (again  by (\ref{rc4}) )} \nonumber\\
&&=(\alpha_V(v)\cdot bc)\cdot\alpha_A^2(a)+(\alpha_V(v)\cdot ca)\cdot\alpha_A^2(b)+
(\alpha_V(v)\cdot ab)\cdot\alpha_A^2(c)\nonumber\\
%\end{eqnarray}
&&\hspace*{-2cm}\mbox{and thus, we get (\ref{rc2}). Finally, (\ref{rc3}) is proved as follows:}\nonumber\\
%\begin{eqnarray}
&&\alpha_V(v\cdot a)\cdot\alpha_A(bc)+\alpha_V(v\cdot b)
\cdot\alpha_A(ca)+\alpha_V(v\cdot c)\cdot\alpha_A(ab)\nonumber\\
&=&\alpha_V(v\cdot a)\cdot\alpha_A(b)\alpha_A(c)+\alpha_V(v\cdot b)
\cdot\alpha_A(c)\alpha_A(a)+\alpha_V(v\cdot c)\cdot\alpha_A(a)\alpha_A(b) (\mbox{ multiplicativity} )\nonumber\\
&=&((v\cdot a)\cdot\alpha_A(b))\cdot\alpha_A^2(c)+((v\cdot a)\cdot\alpha_A(c))\cdot\alpha_A^2(b)+
((v\cdot b)\cdot\alpha_A(c))\cdot\alpha_A^2(a)\nonumber\\
&&+((v\cdot b)\cdot\alpha_A(a))\cdot\alpha_A^2(c)+((v\cdot c)\cdot\alpha_A(a))\cdot\alpha_A^2(b)+((v\cdot c)\cdot\alpha_A(b))\cdot\alpha_A^2(a) \mbox{ (by (\ref{rc4}) )} \nonumber\\
&=&((v\cdot a)\cdot\alpha_A(b))\cdot\alpha_A^2(c)+[\alpha_V(v)\cdot ac-((v\cdot c)\cdot\alpha_A(a)]\cdot\alpha_A^2(b)+
((v\cdot b)\cdot\alpha_A(c))\cdot\alpha_A^2(a)\nonumber\\
&&+((v\cdot b)\cdot\alpha_A(a))\cdot\alpha_A^2(c)+((v\cdot c)\cdot\alpha_A(a))\cdot\alpha_A^2(b)+((v\cdot c)\cdot\alpha_A(b))\cdot\alpha_A^2(a) \mbox{ ( again by (\ref{rc4}) )} \nonumber\\
&=&((v\cdot a)\cdot\alpha_A(b))\cdot\alpha_A^2(c)+\alpha_V^2(v)\cdot((ac)\alpha_A(b)) -(\alpha_V(v)\cdot\alpha_A(b))\cdot\alpha_A(ac)+
((v\cdot b)\cdot\alpha_A(c))\cdot\alpha_A^2(a)\nonumber\\
&&+((v\cdot b)\cdot\alpha_A(a))\cdot\alpha_A^2(c)+((v\cdot c)\cdot\alpha_A(b))\cdot\alpha_A^2(a) \mbox{ ( again by (\ref{rc4}) )} \nonumber\\
&=&((v\cdot a)\cdot\alpha_A(b))\cdot\alpha_A^2(c)+\alpha_V^2(v)\cdot((ac)\alpha_A(b)) -(\alpha_V(v\cdot b)\cdot\alpha_A(ac)+
((v\cdot b)\cdot\alpha_A(c))\cdot\alpha_A^2(a)\nonumber
\end{eqnarray}
\begin{eqnarray}
&&+((v\cdot b)\cdot\alpha_A(a))\cdot\alpha_A^2(c)+((v\cdot c)\cdot\alpha_A(b))\cdot\alpha_A^2(a) \mbox{ ( by  (\ref{rc5})  )} \nonumber\\
&=&((v\cdot a)\cdot\alpha_A(b))\cdot\alpha_A^2(c)+\alpha_V^2(v)\cdot((ac)\alpha_A(b)) -
((v\cdot b)\cdot\alpha_A(a))\cdot\alpha_A^2(c)-((v\cdot b)\cdot\alpha_A(c))\cdot\alpha_A^2(a)\nonumber\\
&&+((v\cdot b)\cdot\alpha_A(c))\cdot\alpha_A^2(a)+((v\cdot b)\cdot\alpha_A(a))\cdot\alpha_A^2(c)+((v\cdot c)\cdot\alpha_A(b))\cdot\alpha_A^2(a) \mbox{ ( by  (\ref{rc4})  )} \nonumber\\
&&=((v\cdot a)\cdot\alpha_A(b))\cdot\alpha_A^2(c)+((v\cdot c)\cdot
\alpha_A(b))\cdot\alpha_A^2(a)+\alpha_V^2(v)\cdot((ac)\alpha_A(b))\nonumber
\end{eqnarray}
which is ((\ref{rc3})).\hfill $\square$\\
\\
 Similarly, the following result can be proved.
\begin{Theorem}\label{thmleft}
Let $(A,\mu,\alpha_A)$ be a Hom-Jordan algebra, $(V,\alpha_V)$ be a Hom-module and $\rho_l: A\otimes V\rightarrow V,$ $(v\otimes a\mapsto a\cdot v),$ be a bilinear map satisfying
\begin{eqnarray}
\alpha_V\circ\rho_l=\rho_l\circ(\alpha_A\otimes\alpha_V) \label{lc5}\\
(ab)\cdot\alpha_V(v)=\alpha_A(a)\cdot(b\cdot v)+\alpha_A(b)\cdot(a\cdot v)\label{lc4}
\end{eqnarray}
for all $(a,b)\in A^2$ and $v\in V.$ Then $(V,\alpha, \rho_l)$ is a left $A$-module called
a left special $A$-module.\hfill $\square$
\end{Theorem}
We know that the plus algebra of any Hom-alternative algebra is a Hom-Jordan algebra.
The next result says that any left (resp. right) module over a Hom-alternative algebra is also 
a left (resp. right) module over its plus Hom-algebra.
\begin{Proposition}
Let $(A,\mu,\alpha_A)$ be a Hom-alternative algebra and $(V,\alpha_V)$ be a Hom-module.
\begin{enumerate}
\item If $(V,\alpha_V)$ is a right $A$-module with the structure maps $\rho_r$ 
then $(V,\alpha_V)$ is a right special $A^+$-module with the same
structure map $\rho_r.$ 
\item If $(V,\alpha_V)$ is a left $A$-module with the structure maps $\rho_l$ 
then $(V,\alpha_V)$ is a left special $A^+$-module with the same
structure map $\rho_l.$ 
\end{enumerate}
\end{Proposition}
{\bf Proof.} It suffices to prove (\ref{rc4}) and (\ref{lc4}). 
\begin{enumerate}
\item If $(V,\alpha_V)$ is a right $A$-module with the structure map $\rho_r$ then for all $(x,y,v)\in A\times A\times V,$
$as_{A,V}(v,x,y)=-as_{A,V}(v,y,x)$ by (\ref{lem2}),  i.e.  $\alpha_V(v)\cdot(xy)+\alpha_V(v)\cdot (yx)=(v\cdot x)\cdot\alpha_A(y)+(v\cdot y)\cdot\alpha_A(x).$ Thus
$\alpha_V(v)\cdot(x\ast y)=\alpha_V(v)\cdot(xy)+\alpha_V(v)\cdot(yx)=(v\cdot x)\cdot\alpha_A(y)+(v\cdot y)\cdot\alpha_A(x).$ Therefore $(V,\alpha_V)$ is a right special $A^+$-module by Theorem \ref{thmright}.
\item If $(V,\alpha_V)$ is a left $A$-module with the structure map $\rho_l$  then for all 
$(x,y,v)\in A\times A\times V,$
$as_{A,V}(x,y,v)=-as_{A,V}(y,x,v)$ by  Remark \ref{rem2} and then $(xy)\cdot\alpha_V(v)+(yx)\cdot\alpha_V(v)=
\alpha_A(x)\cdot(y\cdot v)+\alpha_A(y)\cdot(x\cdot v).$  Thus
$(x\ast y)\cdot\alpha_V(v)=(xy)\cdot\alpha_V(v)+(yx)\cdot\alpha_V(v)=\alpha_A(x)\cdot(y\cdot v)+\alpha_A(y)\cdot(x\cdot v).$  Therefore  $(V,\alpha_V)$ is a left special $A^+$-module 
by Theorem \ref{thmleft}.\hfill $\square$\\
\end{enumerate}
Now, we give the definition of a bimodule over a Hom-Jordan algebra.
\begin{Definition}
Let $(A,\mu,\alpha_A)$ be a Hom-Jordan algebra. \\
(i)An $A$-bimodule is a Hom-module $(V,\alpha_V)$ that comes equipped with a left structure map $\rho_l:A\otimes V\longrightarrow V$ ($\rho_l(a\otimes v)=a\cdot v$) and a right structure map $\rho_r:V\otimes A\longrightarrow V$ ($\rho_r(v\otimes a)=v\cdot a$) 
such that the following conditions:
\begin{eqnarray}
%\alpha_V\circ\rho_l=\rho_l\circ(\alpha_A\otimes\alpha_V)\label{r1}\\
%\alpha_V\circ\rho_r=\rho_l\circ(\alpha_V\otimes\alpha_A)\label{r2}\\
\rho_r\circ\tau_1=\rho_l\label{r3}
\end{eqnarray}
\begin{eqnarray}
&&\alpha_V(v\cdot a)\cdot\alpha_A(bc)+\alpha_V(v\cdot b)
\cdot\alpha_A(ca)+\alpha_V(v\cdot c)\cdot\alpha_A(ab)\nonumber\\
&&=(\alpha_V(v)\cdot bc)\cdot\alpha_A^2(a)+(\alpha_V(v)\cdot ca)\cdot\alpha_A^2(b)+
(\alpha_V(v)\cdot ab)\cdot\alpha_A^2(c)\label{r4}
\end{eqnarray}
\begin{eqnarray}
&&\alpha_V(v\cdot a)\cdot\alpha_A(bc)+\alpha_V(v\cdot b)
\cdot\alpha_A(ca)+\alpha_V(v\cdot c)\cdot\alpha_A(ab)\nonumber\\
&&=((v\cdot a)\cdot\alpha_A(b))\cdot\alpha_A^2(c)+((v\cdot c)\cdot
\alpha_A(b))\cdot\alpha_A^2(a)+((ac)\alpha_A(b))\cdot\alpha_V^2(v)\label{r5}
\end{eqnarray}
hold for all $a,b,c\in A$ and $v\in V.$\\
(ii) A morphism $f :(V, \alpha_V )\longrightarrow (W, \alpha_W )$ of $A$-bimodules is a morphism of the underlying Hom-modules such that
\begin{eqnarray}
f\circ\rho_l=\rho_l\circ(Id_A\otimes f)\nonumber\\
f\circ\rho_r=\rho_l\circ(f \otimes Id_A)\nonumber
\end{eqnarray}
\end{Definition}
 In term of the module Hom-associator,  using the fact that the structure maps are morphisms and the relation (\ref{r3}), the relations (\ref{r4}) and (\ref{r5}) are respectively
\begin{eqnarray}
\circlearrowleft_{(a,b,c)}as_{A,V}(\alpha_A(a),\alpha_V(v),bc)=0\label{nr4}\\
as_{A,V}(v\cdot a,\alpha_A(b),\alpha_A(c))+as_{A,V}(v\cdot c,\alpha_A(b),\alpha_A(a))+
as_{A,V}(ac,\alpha_A(b),\alpha_V(v))=0\label{nr5}
\end{eqnarray}
\begin{Remark}\label{rem1}
(i) Since $\rho_r\circ\tau_1=\rho_l,$ nothing is lost in dropping one of the compositions. Thus the term Hom-Jordan module can be used for Hom-Jordan bimodule.\\
(ii) Since the field is of characteristic $0,$ the identity (\ref{nr4})
implies $as_{A,V}(\alpha_A(a),\alpha_V(v),a^2)=0$
\end{Remark}
\begin{Example} Here are some examples of $A$-modules.\\
(1) Let $(A,\mu,\alpha_A)$ be a Hom-Jordan algebra. Then $(A,\alpha_A)$ is an $A$-bimodule where the structure maps are $\rho_l=\rho_r=\mu.$\\
(2) If $(A,\mu)$ is a Jordan algebra and $M$ is an $A$-bimodule \cite{Jacob2} in the usual sense then $(M,Id_M)$ is an $\mathbb{A}$-bimodule where $\mathbb{A}=(A,\mu, Id_A)$ is a Hom-Jordan algebra.
\end{Example}
The next result shows that a special left and right module over a  Hom-Jordan algebra, has 
a bimodule structure over this Hom-algebra under a specific condition.
\begin{Theorem}\label{ModPlus}
Let $(A,\mu,\alpha_A)$ be a Hom-Jordan algebra and $(V,\alpha_V)$ be both a left  and a right 
special $A$-module with the structure maps $\rho_1$ and $\rho_2$ respectively such that the Hom-associativity (or operator commutativity) condition
\begin{eqnarray}
\rho_2\circ(\rho_1\otimes\alpha_A)=\rho_1\circ(\alpha_A\otimes\rho_2) \label{associat}
\end{eqnarray}
holds. Define the bilinear maps $\rho_l: A\otimes V\longrightarrow V$ and $\rho_r: V\otimes A\longrightarrow V$ by
\begin{eqnarray}
\rho_l=\rho_1+\rho_2\circ\tau_1 \mbox{ and } \rho_r=\rho_1\circ\tau_2+\rho_2 \label{StrucPlus}
\end{eqnarray}
Then $(V, \alpha_V,\rho_l,\rho_r)$ is an $A$-bimodule.
\end{Theorem}
{\bf Proof.} It is clear that $\rho_l$ and $\rho_r$ are structure maps  and (\ref{r3}) holds. To prove relations (\ref{r4}) and (\ref{r5}), let put $\rho_l(a\otimes v):=a\diamond v$ \  i.e. $a\diamond v=a\cdot v+v\cdot a$ for all $(a,v)\in A\times V.$ We have then $\rho_r(v\otimes a):= v\diamond a=a\cdot v+v\cdot a$ 
for all $(a,v)\in A\times V.$ Therefore for all $(a,b,v)\in A\times A\times V,$
\begin{eqnarray}
&&\alpha_V(v\diamond a)\diamond\alpha_A(bc)+\alpha_V(v\diamond b)
\diamond\alpha_A(ca)+\alpha_V(v\diamond c)\diamond\alpha_A(ab)\nonumber\\
&=&\alpha_V(v\cdot a)\cdot\alpha_A(bc)+\alpha_V(a\cdot v)\cdot\alpha_A(bc)+\alpha_A(bc
)\cdot\alpha_V(v\cdot a)+\alpha_A(bc)\cdot\alpha_V(a\cdot v)\nonumber\\
&&+\alpha_V(v\cdot b)\cdot\alpha_A(ca)+\alpha_V(b\cdot v)\cdot\alpha_A(ca)+\alpha_A(ca)\cdot\alpha_V(v\cdot b)
+\alpha_A(ca)\cdot\alpha_V(b\cdot v)\nonumber\\
&&+\alpha_V(v\cdot c)\cdot\alpha_A(ab)+\alpha_V(c\cdot v)\cdot\alpha_A(ab)+\alpha_A(ab)\cdot\alpha_V(v\cdot c)+\alpha_A(ab)\cdot\alpha_V(c\cdot v) \nonumber\\
&&\mbox{ (by a straightforward computation ) }\nonumber\\
&=&\{\alpha_V(v\cdot a)\cdot\alpha_A(bc)+\alpha_V(v\cdot b)\cdot\alpha_A(ca)+\alpha_V(v\cdot c)\cdot\alpha_A(ab)\}+\{\alpha_A(bc)\cdot(\alpha_V(v)\cdot \alpha_A(a))\nonumber\\
&&+\alpha_A(ca)\cdot(\alpha_V(v)\cdot \alpha_A(b))+\alpha_A(ab)\cdot(\alpha_V(v)\cdot \alpha_A(c)) \}
+\{\alpha_A(bc)\cdot\alpha_V(a\cdot v)+\nonumber\\
&&\alpha_A(ca)\cdot\alpha_V(b\cdot v)+\alpha_A(ab)\cdot\alpha_V(c\cdot v) \}
+\{(\alpha_A(a)\cdot \alpha_V(v))\cdot\alpha_A(bc)+(\alpha_A(b)\cdot \alpha_V(v))\cdot\alpha_A(ca)\nonumber\\
&&+(\alpha_A(c)\cdot \alpha_V(v))\cdot\alpha_A(ab)\}\mbox{ (rearranging terms and noting that $\rho_1$ and $\rho_2$ are morphisms  )}\nonumber
\end{eqnarray}
\begin{eqnarray}
&=&\{ (\alpha_V(v)\cdot bc)\cdot\alpha_A^2(a)+(\alpha_V(v)\cdot ca)\cdot\alpha_A^2(b)+
(\alpha_V(v)\cdot ab)\cdot\alpha_A^2(c)\}+\{(bc\cdot\alpha_V(v))\cdot\alpha_A^2(a)\nonumber\\
&&+(ca\cdot\alpha_V(v))\cdot\alpha_A^2(b)+(ab\cdot\alpha_V(v))\cdot\alpha_A^2(c)\}
+\{\alpha_A^2(a)\cdot(bc\cdot\alpha_V(v))+\alpha_A^2(b)\cdot(ca\cdot\alpha_V(v))\nonumber\\
&&+\alpha_A^2(c)\cdot(ab\cdot\alpha_V(v))\}+\{\alpha_A^2(a)\cdot(\alpha_V(v)\cdot bc)+
 \alpha_A^2(b)\cdot(\alpha_V(v)\cdot ca)+\alpha_A^2(c)\cdot(\alpha_V(v)\cdot ab)\}\nonumber\\
&&\mbox{ ( by (\ref{rc2}), (\ref{lc2}) and (\ref{associat}) }\nonumber\\
&=&\{ (\alpha_V(v)\diamond bc)\cdot\alpha_A^2(a)+(\alpha_V(v)\diamond ca)\cdot\alpha_A^2(b)+
(\alpha_V(v)\diamond ab)\cdot\alpha_A^2(c)\}+\{\alpha_A^2(a)\cdot(\alpha_V(v)\diamond bc)\nonumber\\
&&+ \alpha_A^2(b)\cdot(\alpha_V(v)\diamond ca)+\alpha_A^2(c)\cdot(\alpha_V(v)\diamond ab)\} \mbox{ ( by the definition of $\diamond$ )}\nonumber\\
&=& (\alpha_V(v)\diamond bc)\diamond\alpha_A^2(a)+(\alpha_V(v)\diamond ca)\diamond\alpha_A^2(b)+
(\alpha_V(v)\diamond ab)\diamond\alpha_A^2(c)\nonumber\\
&&\mbox{ ( again by the definition of $\diamond$ )}\nonumber\\
%\end{eqnarray}
&&\hspace*{-2cm} \mbox{Therefore we get (\ref{r4}). Finally, we have:}\nonumber\\
%\begin{eqnarray}
&&\alpha_V(v\diamond a)\diamond\alpha_A(bc)+\alpha_V(v\diamond b)
\diamond\alpha_A(ca)+\alpha_V(v\diamond c)\diamond\alpha_A(ab)\nonumber\\
&=&\alpha_V(v\cdot a)\cdot\alpha_A(bc)+\alpha_V(a\cdot v)\cdot\alpha_A(bc)+\alpha_A(bc
)\cdot\alpha_V(v\cdot a)+\alpha_A(bc)\cdot\alpha_V(a\cdot v)\nonumber\\
&&+\alpha_V(v\cdot b)\cdot\alpha_A(ca)+\alpha_V(b\cdot v)\cdot\alpha_A(ca)+\alpha_A(ca)\cdot\alpha_V(v\cdot b)
+\alpha_A(ca)\cdot\alpha_V(b\cdot v)\nonumber\\
&&+\alpha_V(v\cdot c)\cdot\alpha_A(ab)+\alpha_V(c\cdot v)\cdot\alpha_A(ab)+\alpha_A(ab)\cdot\alpha_V(v\cdot c)+\alpha_A(ab)\cdot\alpha_V(c\cdot v)\nonumber\\ 
&&\mbox{ (by a straightforward computation ) }\nonumber\\
&=&\{\alpha_V(v\cdot a)\cdot\alpha_A(bc)+\alpha_V(v\cdot b)\cdot\alpha_A(ca)+\alpha_V(v\cdot c)\cdot\alpha_A(ab)\}+\{(\alpha_V(a\cdot v)\cdot\alpha_A(b)\alpha_A(c)\nonumber\\
&& +(\alpha_V(b\cdot v))\cdot\alpha_A(c)\alpha_A(a)+(\alpha_V(c\cdot v)\cdot\alpha_A(a)\alpha_A(b)\}
+\{\alpha_A(bc)\cdot\alpha_V(a\cdot v)+\alpha_A(ca)\cdot\alpha_V(b\cdot v)\nonumber\\
&&+\alpha_A(ab)\cdot\alpha_V(c\cdot v) \}+\{\alpha_A(b)\alpha_A(c)\cdot\alpha_V(v\cdot a)+\alpha_A(c)\alpha_A(a)\cdot\alpha_V(v\cdot b)\nonumber\\ 
&&+\alpha_A(a)\alpha_A(b)\cdot\alpha_V(v\cdot c)\}\mbox{ (rearranging terms and using the multiplicativity of $\alpha_A$  )}\nonumber\\
&=&\{\underbrace{((v\cdot a)\cdot\alpha_A(b))\cdot\alpha_A^2(c)}_{1}+\underbrace{((v\cdot c)\cdot
\alpha_A(b))\cdot\alpha_A^2(a)}_{2}+\underbrace{\alpha_V^2(v)\cdot((ac)\alpha_A(b))}_{5}\}\nonumber\\
&&+\{\underbrace{((a\cdot v)\cdot\alpha_A(b))\cdot\alpha_A^2(c)}_{1}+((a\cdot v)\cdot\alpha_A(c))\cdot\alpha_A^2(b)+((b\cdot v)\cdot\alpha_A(c))\cdot\alpha_A^2(a)
\nonumber\\
&&+((b\cdot v)\cdot\alpha_A(a))\cdot\alpha_A^2(c)+((c\cdot v)\cdot\alpha_A(a))\cdot\alpha_A^2(b)+\underbrace{((c\cdot v)\cdot\alpha_A(b))\cdot\alpha_A^2(a)}_{2}\}\nonumber\\
&&+\{\underbrace{\alpha_A^2(c)\cdot(\alpha_A(b)\cdot(a\cdot v))}_{3}+\underbrace{\alpha_A^2(a)\cdot(\alpha_A(b)\cdot(c\cdot v))}_{4}+\underbrace{((ac)\alpha_A(b))\cdot\alpha_V^2(v)}_{5}\}\nonumber\\
&&+\{\alpha_A^2(b)\cdot(\alpha_A(c)\cdot(v\cdot a))+\underbrace{\alpha_A^2(c)\cdot(\alpha_A(b)\cdot(v\cdot a))}_{3}+
\alpha_A^2(a)\cdot(\alpha_A(c)\cdot(v\cdot b))\nonumber\\
&&+\alpha_A^2(c)\cdot(\alpha_A(a)\cdot(v\cdot b))+\underbrace{\alpha_A^2(a)\cdot(\alpha_A(b)\cdot(v\cdot c))}_{4}+\alpha_A^2(b)\cdot(\alpha_A(a)\cdot(v\cdot c))\}\nonumber\\
&&\mbox{ ( by (\ref{rc3}), (\ref{lc3}), (\ref{rc4}) and (\ref{lc4}) ) }\nonumber\\
&=&((v\diamond a)\cdot\alpha_A(b))\cdot\alpha_A^2(c)+((v\diamond c)\cdot
\alpha_A(b))\cdot\alpha_A^2(a)+\alpha_A^2(c)\cdot(\alpha_A(b)\cdot(v\diamond a ))\nonumber\\
&&+\alpha_A^2(a)\cdot(\alpha_A(b)\cdot(v\diamond c))+\alpha_V^2(v)\diamond((ac)\alpha_A(b))+((a\cdot v)\cdot\alpha_A(c))\cdot\alpha_A^2(b)\nonumber\\
&&+((b\cdot v)\cdot\alpha_A(c))\cdot\alpha_A^2(a)+((b\cdot v)\cdot\alpha_A(a))\cdot\alpha_A^2(c)+((c\cdot v)\cdot\alpha_A(a))\cdot\alpha_A^2(b)\nonumber\\
&&+\alpha_A^2(b)\cdot(\alpha_A(c)\cdot(v\cdot a))+\alpha_A^2(a)\cdot(\alpha_A(c)\cdot(v\cdot b))+\alpha_A^2(c)\cdot(\alpha_A(a)\cdot(v\cdot b))\nonumber\\
&&+\alpha_A^2(b)\cdot(\alpha_A(a)\cdot(v\cdot c))\nonumber\\
&=&((v\diamond a)\cdot\alpha_A(b))\cdot\alpha_A^2(c)+((v\diamond c)\cdot
\alpha_A(b))\cdot\alpha_A^2(a)+\alpha_A^2(c)\cdot(\alpha_A(b)\cdot(v\diamond a ))\nonumber\\
&&+\alpha_A^2(a)\cdot(\alpha_A(b)\cdot(v\diamond c))+\alpha_V^2(v)\diamond((ac)\alpha_A(b))+(\alpha_A(a)\cdot( v\cdot c))\cdot\alpha_A^2(b)\nonumber\\
&&+(\alpha_A(b)\cdot( v\cdot c))\cdot\alpha_A^2(a)+(\alpha_A(b)\cdot (v\cdot a))\cdot\alpha_A^2(c)+(\alpha_A(c)\cdot (v\cdot a))\cdot\alpha_A^2(b)\nonumber\\
&&+\alpha_A^2(b)\cdot((c\cdot v)\cdot \alpha_A(a))+
\alpha_A^2(a)\cdot((c\cdot v)\cdot \alpha_A(b))+\alpha_A^2(c)\cdot((a\cdot v)\cdot \alpha_A(b))\nonumber\\
&&+\alpha_A^2(b)\cdot((a\cdot v)\cdot \alpha_A(c)) \mbox{ ( by  (\ref{associat}) )}\nonumber\\
&=&((v\diamond a)\cdot\alpha_A(b))\cdot\alpha_A^2(c)+((v\diamond c)\cdot
\alpha_A(b))\cdot\alpha_A^2(a)+\alpha_A^2(c)\cdot(\alpha_A(b)\cdot(v\diamond a ))\nonumber
\end{eqnarray}
\begin{eqnarray}
&&+\alpha_A^2(a)\cdot(\alpha_A(b)\cdot(v\diamond c))+\alpha_V^2(v)\diamond((ac)\alpha_A(b))+\underbrace{\alpha_A^2(a)\cdot( (v\cdot c)\cdot\alpha_A(b))}_{6}\nonumber\\
&&+\underbrace{(\alpha_A(b)\cdot( v\cdot c))\cdot\alpha_A^2(a)}_{7}+\underbrace{(\alpha_A(b)\cdot (v\cdot a))\cdot\alpha_A^2(c)}_{8}+\underbrace{\alpha_A^2(c)\cdot( (v\cdot a)\cdot\alpha_A(b))}_{9}
\nonumber\\
&&+\underbrace{(\alpha_A(b)\cdot(c\cdot v))\cdot \alpha_A^2(a)}_{7}+
\underbrace{\alpha_A^2(a)\cdot((c\cdot v)\cdot \alpha_A(b))}_{6}+\underbrace{\alpha_A^2(c)\cdot((a\cdot v)\cdot \alpha_A(b))}_{9}\nonumber\\
&&+\underbrace{(\alpha_A(b)\cdot(a\cdot v))\cdot \alpha_A^2(c))}_ {8} \mbox{ ( again by  (\ref{associat}) )}\nonumber\\
&=&\underbrace{((v\diamond a)\cdot\alpha_A(b))\cdot\alpha_A^2(c)}_ {10}+\underbrace{((v\diamond c)\cdot
\alpha_A(b))\cdot\alpha_A^2(a)}_{11}+\underbrace{\alpha_A^2(c)\cdot(\alpha_A(b)\cdot(v\diamond a ))}_ {13}\nonumber\\
&&+\underbrace{\alpha_A^2(a)\cdot(\alpha_A(b)\cdot(v\diamond c))}_{12}+\alpha_V^2(v)\diamond((ac)\alpha_A(b))+\underbrace{\alpha_A^2(a)\cdot( (v\diamond c)\cdot\alpha_A(b))}_ {12}\nonumber\\
&&+\underbrace{(\alpha_A(b)\cdot( v\diamond c))\cdot\alpha_A^2(a)}_{11}+
\underbrace{(\alpha_A(b)\cdot (v\diamond a))\cdot\alpha_A^2(c)}_ {10}+
\underbrace{\alpha_A^2(c)\cdot( (v\diamond a)\cdot\alpha_A(b))}_{13}\nonumber\\
&=&((v\diamond a)\diamond\alpha_A(b))\cdot\alpha_A^2(c)+((v\diamond c)\diamond
\alpha_A(b))\cdot\alpha_A^2(a)
+\alpha_A^2(a)\cdot( (v\diamond c)\diamond\alpha_A(b))
\nonumber\\
&&+\alpha_A^2(c)\cdot( (v\diamond a)\diamond\alpha_A(b))+\alpha_V^2(v)\diamond((ac)\alpha_A(b))\nonumber\\
&=&((v\diamond a)\diamond\alpha_A(b))\diamond\alpha_A^2(c)+((v\diamond c)\diamond
\alpha_A(b))\diamond\alpha_A^2(a)
+\alpha_V^2(v)\diamond((ac)\alpha_A(b))\nonumber
\end{eqnarray}
which is (\ref{r5}).\hfill $\square$\\

The following result will be used below. It gives a relation between modules over Hom-alternative 
algebras and special modules over Hom-Jordan algebras.
\begin{Lemma}\label{HomAssModule}
Let $(A,\mu,\alpha_A)$ be a Hom-associative algebra and $(V,\alpha_V)$ be a Hom-module.
\begin{enumerate}
\item If $(V,\alpha_V)$ is a right $A$-module with the structure maps $\rho_r$ 
then $(V,\alpha_V)$ is a right special $A^+$-module with the same
structure map $\rho_r.$ 
\item If $(V,\alpha_V)$ is a left $A$-module with the structure maps $\rho_l$ 
then $(V,\alpha_V)$ is a left special $A^+$-module with the same
structure map $\rho_l.$ 
\end{enumerate}
\end{Lemma}
{\bf Proof.} It also suffices to prove (\ref{rc4}) and (\ref{lc4}). 
\begin{enumerate}
\item If $(V,\alpha_V)$ is a right $A$-module with the structure map $\rho_r$ then for all $(x,y,v)\in A\times A\times V,$
$\alpha_V(v)\cdot(a\ast b)=\alpha_V(v)\cdot(ab)+\alpha_V(v)\cdot(ab)=(v\cdot a)\cdot\alpha_A(b)+(v\cdot b)\cdot\alpha_A(a)$
where the last equality holds by (\ref{RightAssMod}). Then $(V,\alpha_V)$ is a right special $A^+$-module.
\item If $(V,\alpha_V)$ is a left $A$-module with the structure map $\rho_l$  then for all 
$(x,y,v)\in A\times A\times V,$ 
$(a\ast b)\cdot\alpha_V(v)=(ab)\cdot\alpha_V(v)+(ba)\cdot\alpha_V(v)=\alpha_A(a)\cdot(b\cdot v)+\alpha_A(b)\cdot(a\cdot v)$ where the last equality holds by (\ref{LeftAssMod}). Then  $(V,\alpha_V)$ is a left special $A^+$-module.\hfill $\square$
\end{enumerate}
Now we prove that a module over a Hom-associative algebra gives rise to a module over its plus Hom-algebra.
\begin{Proposition}\label{bhamHJ}
Let $(A,\mu,\alpha_A)$ be a Hom-associative algebra and $(V, \rho_1, \rho_2,\alpha_V)$  be an 
$A$-bimodule. Then $(V, \rho_l, \rho_r,\alpha_V)$  is an $A^+$-bimodule where $\rho_l$ and $\rho_r$ are defined as in (\ref{StrucPlus}).
\end{Proposition}
{\bf Proof.} The proof follows from Lemma \ref{HomAssModule} , the Hom-associativity condition 
(\ref{Assoc1}) and  Theorem \ref{ModPlus}.\hfill $\square$\\

The following elementary but important result will be used below. It gives a property of a module 
Hom-associator.
\begin{Lemma}
Let $(A,\mu,\alpha_A)$ be a Hom-Jordan algebra and  $(V,\alpha_V)$  be an $A$-bimodule with the structure maps $\rho_l$ and $\rho_r$. Then 
\begin{eqnarray}
\alpha_V^n\circ as_{A,V}\circ Id_{A\otimes V\otimes A}=as_{A,V}\circ(\alpha_A^{\otimes n}\otimes\alpha_V^{\otimes n}\otimes\alpha_A^{\otimes n})\label{lem1}
\end{eqnarray}
\end{Lemma}
{\bf Proof.} Using twice the fact that $\rho_l$ and $\rho_r$ are morphisms of Hom-modules, we get
\begin{eqnarray}
\alpha_V^n\circ as_{A,V}\circ Id_{A\otimes V\otimes A}&=&\alpha_V^n\circ(\rho_r\circ
(\rho_l\otimes\alpha_A)-\rho_l\circ(\alpha_A\otimes\rho_r))\nonumber\\
 &=&\alpha_V^n\circ\rho_r\circ(\rho_l\otimes\alpha_A)
 -\alpha_V^n\circ\rho_l\circ(\alpha_A\otimes\rho_r) \mbox{ (linearity of $\alpha_V^n$)}\nonumber\\
 &=&\rho_r\circ(\alpha_V^n\circ\rho_l\otimes\alpha_A^{n+1})-
 \rho_l\circ(\alpha_A^{n+1}\otimes\alpha_V^n\circ\rho_r) %\mbox{ (by (\ref{r1}) and (\ref{r2}) )}
 \nonumber\\
 &=&\rho_r\circ(\rho_l\circ(\alpha_A^n\otimes\alpha_V^n)\otimes\alpha_A^{n+1})-
 \rho_l\circ(\alpha_A^{n+1}\otimes\rho_r\circ(\alpha_V^n\otimes \alpha_A^n))\nonumber\\
   %&& \mbox{ ( again by (\ref{r1}) and (\ref{r2}) )}\nonumber\\
   &=&(\rho_r\circ(\rho_l\otimes\alpha_A)-\rho_l\circ(\alpha_A\otimes\rho_r))\circ
   (\alpha_A^{\otimes n}\otimes\alpha_V^{\otimes n}\otimes\alpha_A^{\otimes n})\nonumber\\
   &=&as_{A,V}\circ(\alpha_A^{\otimes n}\otimes\alpha_V^{\otimes n}\otimes\alpha_A^{\otimes n})\nonumber
\end{eqnarray}
That ends the proof.\hfill $\square$\\
\\
The next result is similar to Proposition \ref{sma} . Contrary to Proposition \ref{sma}, 
 an additional condition is needed.
\begin{Proposition}\label{HJB-HJB}
Let $(A,\mu,\alpha_A)$ be a Hom-Jordan algebra and  $(V,\alpha_V)$  be an $A$-bimodule with the structure maps $\rho_l$ and $\rho_r$. Suppose that there exists $n\in\mathbb{N}$ such that $\alpha_V^{n}=Id_V.$ Then the maps
\begin{eqnarray}
\rho_l^{(n)}=\rho_l\circ(\alpha_A^n\otimes Id_V)\label{nmj1}\\
\rho_r^{(n)}=\rho_r\circ(Id_V\otimes \alpha_A^n)\label{nmj2}
\end{eqnarray} 
give the Hom-module $(V,\alpha_V)$ the structure of an $A$-bimodule that we denote 
by $V^{(n)}$
\end{Proposition}
{\bf Proof.} Since the structure map $\rho_l$ is a morphism of Hom-modules, we get:
\begin{eqnarray}
\alpha_V\circ\rho_l^{(n)}&=&\alpha_V\circ\rho_l\circ(\alpha_A^n\otimes Id_V) 
\mbox{ ( by (\ref{nmj1}) )}\nonumber\\
&=&\rho_l\circ(\alpha_A^{n+1}\otimes\alpha_V)\\ % \mbox{ ( by (\ref{r1}) )}\nonumber
&=&\rho_l\circ(\alpha_A^n\otimes Id_V)\circ(\alpha_A\otimes\alpha_V)\nonumber\\
&=&\rho_l^{(n)}\circ(\alpha_A\otimes\alpha_V)\nonumber
\end{eqnarray}
Then, $\rho_l^{(n)}$ is a morphism. Similarly we get that $\rho_r^{(n)}$ is a morphism and  that (\ref{r3}) holds for $V^{(n)}.$ Next, we compute
\begin{eqnarray}
&&\circlearrowleft_{(a,b,c)}as_{A,V^{(n)}}(\alpha_A(a),\alpha_V(v),ab)\nonumber\\
&&=\circlearrowleft_{(a,b,c)}\{\rho_r^{(n)}(\rho_l^{(n)}(\alpha_A(a),\alpha_V(v)),\alpha_A(bc))
-\rho_l^{(n)}(\alpha_A^2(a),\rho_r^{(n)}(\alpha_V(v),bc))\}\nonumber\\
&&=\circlearrowleft_{(a,b,c)}\{\rho_r(\rho_l^{(n)}(\alpha_A(a),\alpha_V(v)),\alpha_A^{n+1}(bc))-
\rho_l(\alpha_A^{n+2}(a),\rho_r^{(n)}(\alpha_V(v),bc))\}\nonumber\\
&&=\circlearrowleft_{(a,b,c)}\{\rho_r(\rho_l(\alpha_A^{n+1}(a),\alpha_V(v)),\alpha_A^{n+1}(bc))-
\rho_l(\alpha_A^{n+2}(a),\rho_r(\alpha_V(v),\alpha_A^n(bc))\}\nonumber\\
&&=\circlearrowleft_{(a,b,c)}\{\rho_r(\rho_l(\alpha_A^{n+1}(a),\alpha_V(v)),\alpha_A(\alpha_A^n(bc)))-
\rho_l(\alpha_A(\alpha_A^{n+1}(a)),\rho_r(\alpha_V(v),\alpha_A^n(bc))\}\nonumber\\
&&=\circlearrowleft_{(a,b,c)}as_{A,V}(\alpha_A^{n+1}(a),\alpha_V(v),\alpha_A^n(bc))\nonumber\\
&&=\circlearrowleft_{(a,b,c)}as_{A,V}(\alpha_A^{n+1}(a),\alpha_V^{n+1}(v),\alpha_A^n(bc)) \mbox{ ( by the hypothesis $\alpha_V=\alpha_V^{n+1}$ )}\nonumber\\
&&= \alpha_V^n(\circlearrowleft_{(a,b,c)}as_{A,V}(\alpha_A(a),\alpha_V(v),bc))
\mbox{ (by (\ref{lem1}) and the linearity of $\alpha_V^n $)} \nonumber\\
&&=0 \mbox{ ( by (\ref{nr4}) in $V )$}\nonumber
\end{eqnarray}
Then we get (\ref{nr4}) for $V^{(n)}.$ Finally remarking that
\begin{eqnarray}
&&as_{A,V^{(n)}}(\rho_r^n(v,a),\alpha_A(b),\alpha_A(c))\nonumber\\
&&=as_{A,V^{(n)}}(v\cdot\alpha_A^n(a),\alpha_A(b),\alpha_A(c))\nonumber\\
&&=\rho_r^n(\rho_r^n(v\cdot\alpha_A^n(a),\alpha_A(b),\alpha_A^2(c))-
\rho_r^n(\alpha_V(v)\cdot\alpha_A^{n+1}(a),\mu(\alpha_A(b),\alpha(_A(c))\nonumber\\
&&=\rho_r(\rho_r(v\cdot\alpha_A^n(a),\alpha_A^{n+1}(b),\alpha_A^{n+2}(c))-
\rho_r(\alpha_V(v)\cdot\alpha^{n+1}(a),\mu(\alpha_A^{n+1}(b),\alpha_A^{n+1}(c))
\nonumber\\
&&=\alpha_{A,V}(v\cdot\alpha_A^n(a),\alpha_A^{n+1}(b),\alpha_A^{n+1}(c))\nonumber
\end{eqnarray} 
and similary
\begin{eqnarray}
as_{A,V^{(n)}}(\rho_r^n(v,c),\alpha_A(b),\alpha_A(a))
&=&as_{A,V}(v\cdot\alpha_A^n(c),\alpha_A^{n+1}(b),\alpha_A^{n+1}(a))\nonumber\\
as_{A,V^{(n)}}(ac,\alpha_A(b),\alpha_V(v))
&=&as_{A,V}(\alpha_A^n(a)\alpha_A^n(c),\alpha_A^{n+1}(b),\alpha_V(v))\nonumber
\end{eqnarray}
(\ref{nr5}) is proved for $V^{(n)}$ as it follows:
\begin{eqnarray}
&&as_{A,V^{(n)}}(\rho_r^n(v,a),\alpha_A(b),\alpha_A(c))+
as_{A,V^{(n)}}(\rho_r^n(v,c),\alpha_A(b),\alpha_A(a))+
as_{A,V^{(n)}}(ac,\alpha_A(b),\alpha_V(v))\nonumber\\
&&=\alpha_V(v\cdot\alpha_A^n(a),\alpha_A^{n+1}(b),\alpha_A^{n+1}(c))
+as_{A,V}(v\cdot\alpha_A^n(c),\alpha_A^{n+1}(b),\alpha_A^{n+1}(a))\nonumber\\
&&+as_{A,V}(\alpha_A^n(a)\alpha_A^n(c),\alpha_A^{n+1}(b),\alpha_V(v))\nonumber\\
&&=\alpha_V(v\cdot\alpha_A^n(a),\alpha_A(\alpha_A^n(b)),\alpha_A(\alpha_A^n(c)))
+as_{A,V}(v\cdot\alpha_A^n(c),\alpha_A(\alpha_A^n(b)),\alpha_A(\alpha_A^n(a)))\nonumber\\
&&+as_{A,V}(\alpha_A^n(a)\alpha_A^n(c),\alpha_A(\alpha_A^n(b)),\alpha_V(v))\nonumber\\
&&=0 \mbox{ ( by (\ref{nr5}) in $V$) }\nonumber
\end{eqnarray}
We conclude that $V^{(n)}$ is an $A$-bimodule. \hfill $\square$\\
\begin{Corollary}
Let $(A,\mu,\alpha_A)$ be a Hom-Jordan algebra and  $(V,\alpha_V)$  be an $A$-bimodule with the structure maps $\rho_l$ and $\rho_r$ such that $\alpha_V$ is an involution. Then $(V,\alpha_V)$ is an $A$-bimodule with the structure maps
$\rho_l^{(2)}=\rho_l\circ(\alpha_A^2\otimes Id_V)$ and 
$\rho_r^{(2)}=\rho_r\circ(Id_V\otimes\alpha_A^2).$ \hfill $\square$
\end{Corollary}
The following result is similar to theorem \ref{mamHa}. It says that bimodules over Jordan algebras can  be deformed into bimodules over Hom-Jordan algebras via an endomorphism.
\begin{Theorem}\label{HJB-JB}
Let $(A,\mu)$ be a Jordan algebra, $V$ be an $A$-bimodule with the structure maps
$\rho_l$ and $\rho_r$, $\alpha_A$ be an endomorphism of the Jordan algebra $A$ 
and $\alpha_V$ be a linear self-map of $V$ such that $\alpha_V\circ\rho_l=\rho_l\circ(\alpha_A\otimes\alpha_V)$ and 
$\alpha_V\circ\rho_r=\rho_r\circ(\alpha_V\otimes\alpha_A)$\\
Write $A_{\alpha_A}$ for the Hom-Jordan algebra $(A,\mu_{\alpha_A},\alpha_A)$ and 
$V_{\alpha_V}$ for the Hom-module $(V,\alpha_V).$ Then the maps:
\begin{eqnarray}
\tilde{\rho_l}=\alpha_V\circ\rho_l \mbox{ and }
\tilde{\rho_r}=\alpha_V\circ\rho_r
\end{eqnarray}
give the Hom-module $V_{\alpha_V}$ the structure of an $A_{\alpha_A}$-bimodule.
\end{Theorem}
{\bf Proof.} First the fact that $\tilde{\rho_l}$ and $\tilde{\rho_r}$ are morphisms and the relation (\ref{r3}) for $V_{\alpha_V}$ are a straightforward computation. Secondly, remarking that
\begin{eqnarray}
as_{A,V_{\alpha_V}}=\alpha_V^2\circ as_{A,V}\label{ias}
\end{eqnarray}
we first compute
\begin{eqnarray}
&&\circlearrowleft_{(a,b,c)}as_{A,V_{\alpha_V}}(\alpha_A(a),\alpha_V(v),\mu_{\alpha_A}(b,c))\nonumber\\
&=&\circlearrowleft_{(a,b,c)}\alpha_V^2(as_{A,V}(\alpha_A(a),\alpha_V(v),\alpha_A(bc))) \mbox{ (by (\ref{ias}) )}\nonumber\\
&=&\circlearrowleft_{(a,b,c)}\alpha_V^3((as_{A,V}(a,v,bc))\mbox{ (by (\ref{lem1}) )}\nonumber\\
&=&\alpha_V^3(\circlearrowleft_{(a,b,c)}(as_{A,V}(a,v,bc))\nonumber\\
&=& 0  \mbox{( by (\ref{nr4}) in $V$ )}\nonumber\\
%\end{eqnarray}
%\begin{eqnarray}
&&\hspace*{-2,5cm}\mbox{and then, we get (\ref{nr4}) for $V_{\alpha_V}$. Finally  we get}\nonumber\\
&&as_{A,V_{\alpha_V}}(\tilde{\rho_r}(v,a),\alpha_A(b),\alpha_A(c))+as_{A,V_{\alpha_V}}(\tilde{\rho_r}(v,c),\alpha_A(b),\alpha_A(a))\nonumber\\
&&+as_{A,V_{\alpha_V}}(\mu_{\alpha_A}(a,c),\alpha_A(b),\alpha_V(v))\nonumber\\
&&=\alpha_V^2(as_{A,V}(\tilde{\rho_r}(v,a),\alpha_A(b),\alpha_A(c)))+\alpha_V^2(as_{A,V}(\tilde{\rho_r}(v,c),\alpha_A(b),\alpha_A(a)))\nonumber\\
&&+\alpha_V^2(as_{A,V}(\mu_{\alpha_A}(a,c),\alpha_A(b),\alpha_V(v))) \mbox{ (by (\ref{ias}) )}\nonumber\\
&&=\alpha_V^2(as_{A,V}(\alpha_V(v\cdot a),\alpha_A(b),\alpha_A(c)))+\alpha_V^2(as_{A,V}(\alpha_V(v\cdot c),\alpha_A(b),\alpha_A(a)))\nonumber\\
&&+\alpha_V^2(as_{A,V}(\alpha_A(ac),\alpha_A(b),\alpha_V(v))) \nonumber\\
&&=\alpha_V^3(as_{A,V}(v\cdot a,b,c))+\alpha_V^3(as_{A,V}(v\cdot c,b,a))
+\alpha_V^3(as_{A,V}(ac,b,v)) \mbox{ (by (\ref{lem1}) )}\nonumber\\
&&=\alpha_V^3(as_{A,V}(v\cdot a,b,c)+as_{A,V}(v\cdot c,b,a)+as_{A,V}(ac,b,v)) \nonumber\\
&&=0 \mbox{( by (\ref{nr5}) in $V$ )} \mbox{ which is (\ref{nr5}) for $V_{\alpha_V}$. }\nonumber
\end{eqnarray}
Therefore the Hom-module $V_{\alpha_V}$ has an $A_{\alpha_A}$-bimodule structure.\hfill $\square$
\begin{Corollary}
Let $(A,\mu)$ be a Jordan algebra, $V$ be a $A$-bimodule with the structure maps
$\rho_l$ and $\rho_r$, $\alpha_A$ be  an endomorphism of the Jordan algebra $A$ 
and $\alpha_V$ be a linear self-map of $V$ such that $\alpha_V\circ\rho_l=\rho_l\circ(\alpha_A\otimes\alpha_V)$ and 
$\alpha_V\circ\rho_r=\rho_r\circ(\alpha_V\otimes\alpha_A)$\\
Moreover suppose that there exists $n\in\mathbb{N}$ such that $\alpha_V^{n}=Id_V.$ 
Write $A_{\alpha_A}$ for the Hom-Jordan algebra $(A,\mu_{\alpha_A},\alpha_A)$ and 
$V_{\alpha_V}$ for the Hom-module $(V,\alpha_V).$ 
Then the maps:
\begin{eqnarray}
\tilde{\rho_l^n}=\rho_l\circ(\alpha_A^{n+1}\otimes\alpha_V) \mbox{ and }
\tilde{\rho_r^n}=\rho_r\circ(\alpha_V\otimes\alpha_A^{n+1})
\end{eqnarray}
give the Hom-module $V_{\alpha}$ the structure of an $A_{\alpha_A}$-bimodule for 
each $n\in\mathbb{N}$.
\end{Corollary}
{\bf Proof.} The proof follows from Proposition \ref{HJB-HJB} and Theorem \ref{HJB-JB}.
\hfill $\square$\\
\\
Similarly to Hom-alternative algebras, the split null extension determined by the
given bimodule over a Hom-Jordan algebra, is constructed as follows:
\begin{Theorem}
Let $(A,\mu,\alpha_A)$ be a Hom-Jordan algebra and  $(V,\alpha_V)$  be an $A$-bimodule with the structure maps $\rho_l$ and $\rho_r$. If define on $A\oplus V$ the linear maps
$\tilde{\mu}: (A\oplus V)^{\otimes 2}\longrightarrow A\oplus V,$ 
$\tilde{\mu}(a+m,b+n):=ab+a\cdot n+m\cdot b$  and 
$\tilde{\alpha}: A\oplus V\longrightarrow A\oplus V,$ $\tilde{\alpha}(a+m):=\alpha_A(a)+\alpha_V(m)$ then $(A\oplus V,\tilde{\mu},\tilde{\alpha})$ is 
a Hom-Jordan algebra.
\end{Theorem}
{\bf Proof.} First, the commutativity of $\tilde{\mu}$ follows from the one of $\mu.$
Next, the multiplicativity of $\tilde{\alpha}$ with respect to $\tilde{\mu}$ follows from the one of $\alpha$ with respect to $\mu$ and the that that $\rho_l$ and $\rho_r$ are morphisms of Hom-modules.
Finally, we prove the Hom-Jordan identity (\ref{idJord}) for $E=A\oplus V$ as it follows
\begin{eqnarray}
&&as_E(\tilde{\mu}(x+m,x+m),\tilde{\alpha}(y+n),\tilde{\alpha}(x+m))\nonumber\\
&&=\tilde{\mu}(\tilde{ \mu}(\tilde{\mu}(x+m,x+m),\tilde{\alpha}(y+n)),\tilde{\alpha}^2(x+m))
-\tilde{\mu}(\tilde{\alpha}(\tilde{\mu}(x+m,x+m)),\tilde{\mu}(\tilde{\alpha}(y+n),
\tilde{\alpha}(x+m)))\nonumber\\
&&=\tilde{\mu}(\tilde{\mu}(x^2+x\cdot m+m\cdot x,\alpha_A(y)+\alpha_V(n)),\alpha_A^2(x)+\alpha_V^2(m))\nonumber\\
&&-\tilde{\mu}(\alpha_A(x^2)+\alpha_V(x\cdot m)+\alpha_V(m\cdot x),\tilde{\mu}(\alpha_A(y)+\alpha_V(n),\alpha_A(x)+\alpha_V(m)))\nonumber
\end{eqnarray}
\begin{eqnarray}
&&=\tilde{\mu}(x^2\alpha_A(y)+x^2\cdot\alpha_V(n)+(x\cdot m)\cdot\alpha_A(y)+
(m\cdot x)\cdot\alpha_A(y),\alpha_A^2(x)+\alpha_V^2(m))\nonumber\\
&&-\tilde{\mu}(\alpha_A^2(x^2)+\alpha_V(x\cdot m)+\alpha_V(m\cdot x),\alpha_A(y)
\alpha_A(x)+\alpha_A(y)\cdot\alpha_V(m)+\alpha_V(n)\cdot\alpha_A(x))\nonumber\\
&&=(x^2\alpha_A(y))\alpha_A^2(x)+(x^2\alpha_A(y))\cdot\alpha_V^2(m)+(x^2\cdot
\alpha_V(n))\cdot\alpha_A^2(x)+((x\cdot m)\cdot\alpha_A(y))\cdot\alpha_A^2(x)\nonumber\\
&&+((m\cdot x)\cdot\alpha_A(y))\cdot\alpha_A^2(x))-\alpha_A(x^2)(\alpha_A(y)\alpha_A(x))-\alpha_A(x^2)\cdot(\alpha_A(y)\cdot\alpha_V(m))\nonumber\\
&&-\alpha_A(x^2)\cdot(\alpha_V(n)\cdot\alpha_A(x))-
\alpha_V(x\cdot m)\cdot(\alpha_A(y)\alpha_A(x))-
\alpha_V(m\cdot x)\cdot(\alpha_A(y)\alpha_A(x))\nonumber\\
&&=as_A(x^2,\alpha_A(y),\alpha_A(x))+as_{A,V}(x^2,\alpha_A(y),\alpha_V(m))+
as_{A,V}(x^2,\alpha_V(n),\alpha_A(x))+\nonumber\\
&&+as_{A,V}(x\cdot m,\alpha_A(y),\alpha_A(x))+as_{A,V}(m\cdot x,\alpha_A(y),\alpha_A(x))\nonumber\\
&&=\underbrace{as_{A,V}(m\cdot x,\alpha_A(y),\alpha_A(x))+as_{A,V}(m\cdot x,\alpha_A(y),\alpha_A(x))
+as_{A,V}(x^2,\alpha_A(y),\alpha_V(m))}_{0}+\nonumber\\
&&+\underbrace{as_{A,V}(x^2,\alpha_V(n),\alpha_A(x))}_{0}+\underbrace{as_A(x^2,\alpha_A(y),\alpha_A(x))}_{0}\nonumber\\
&&=0, \mbox{ where the first $0$ follows from (\ref{nr5}), the second from (\ref{nr4}) (see Remark \ref{rem1}) and the last}\nonumber
\end{eqnarray}
 from the Hom-Jordan identity (\ref{idJord}) in $A.$ We conclude then that $(A\oplus V,\tilde{\mu},\tilde{\alpha})$ is a Hom-Jordan\\ algebra.\hfill $\square$


\begin{thebibliography}{99}
\bibitem{HT1} H. Atagema, A. Makhlouf and S. D. Silvestrov, {\}it Generalization of n-ary Nambu algebras and beyond,} J. Math.Phys. ,{\bf 50}(8)(2009), 083501.
\bibitem{ibbm} I. Bakayoko and B. Manga {\it Hom-alternative modules and Hom-Poisson comodules,} arXiv:1411.7957v1. 
\bibitem{se} S. Eilenberg, {\it Extensions of general algebras,} Annales de la Societe Polonaise de
Mathematique, {\bf 21} (1948), 125-34.
\bibitem{fgcht}  F. G\"ursey and C.-H. Tze, {\it On The Role of Division, Jordan and Related
Algebras in Particle Physics,} World Scientific, Singapore, 1996.
\bibitem{HAR1} J. T. Hartwig, D. Larsson and S. D. Silvestrov, {\it Deformations of Lie algebras using $\sigma-$derivations,} J. Algebras, \textbf{292} (2006), 314-361. 
\bibitem{Jacob1} N. Jacobson, {\it General representation theorie of Jordan algebras,} Trans. Amer. Math. Soc., \textbf{70} (1951), 509-530.
\bibitem{Jacob2}  N. Jacobson, {\it Structure of alternative and Jordan bimodules,} 
Osaka J. Math. {\bf 6} (1954),1-71.
\bibitem{pjjnew} P. Jordan, J. von Neumann, and E. Wigner, {\it On an algebraic generalization
of the quantum mechanical formalism,} Ann. Math., {\bf 35} (1934), 29-64.
 \bibitem{dlsds1} D. Larsson and S. D. Silvestrov, {\it Quasi-Hom-Lie algebras, Central Extensions
and 2-cocycle-like identities,} J. of Algebra, {\bf 288} (2005), 321-344.
\bibitem{dlsds2} D. Larsson and S. D. Silvestrov, {\it Quasi-Lie algebras,} in "Noncommutative
Geometry and Representation Theory in Mathematical Physics", Contemp.
Math., 391, Amer. Math. Soc., Providence, RI (2005), 241-248.
\bibitem{dlsds3} D. Larsson and S. D. Silvestrov, Quasi-deformations of $sl_2(F)$ using twisted
derivations, Comm. Algebra, {\bf 35} (2007), 4303 – 4318.
\bibitem{MAK0} Makhlouf, A., Silvestrov, S. {\it: Notes on formal deformations of Hom-associative and Hom-Lie algebras.} Forum Math., {\bf 22}(4), 715-739 (2010)
\bibitem{MAK1} A. Makhlouf , {\it Hom-Alternative algebras and Hom-Jordan algebras,} Int. Elect. J. Alg., \textbf{8} (2010), 177-190.
%\bibitem{MAK2} A. Makhlouf, {\it Paradigm of nonassociative Hom-algebras and Hom-%superalgebras,} Proceedings of Jordan structures in algebra and Analysis
%Meeting, 143-177, Editional Circulo Rojo, Almeria, 2010.
\bibitem{MAK3} A. Makhlouf, Silvestrov S.D., {\it Hom-algebra structures,} J. Gen. Lie Theory Appl. \textbf{2} (2008), 51-64.
\bibitem{so} S. Okubo, {\it Introduction to octonion and other non-associative algebras in
physics,} Cambridge Univ. Press, Cambridge, UK, 1995.
\bibitem{Schaf1} R. D. Schafer, {\it Representations of alternative algebras,} Trans. Amer. Math. Soc., \textbf{72}(1952), 1-17.
\bibitem{YS} Yunhe Sheng, {\it Representation of hom-Lie algebras.} Algebr. Represent. Theory {\bf 15} (2012), no. 6, 1081-1098.
\bibitem{tasfdv} T.A. Springer and F.D. Veldkamp, {\it Octonions, Jordan Algebras, and 
Exceptional Groups,} Springer, Berlin, 2000.
\bibitem{jtrmw} J. Tits and R.M. Weiss, {\it Moufang Polygons,} Springer-Verlag, Berlin, 2002.
\bibitem{YAU4} D. Yau, {\it Hom-algebras as deformations and homology,} arXiv:0712.3515v1.
\bibitem{YAU2}  D. Yau, Module Hom-algebras, arXiv:0812.4695v1.
\bibitem{YAU3} D. Yau, {\it Hom-Maltsev, Hom-alternative and Hom-Jordan algebras,} International Electronic Journal of algebras, \textbf{11} (2012), 
177-217. 

%%%%%%%%%%%%%%%%%%%%%%%%%%%%%%%%%%%%%%%%%%%%%%%%%%%%%%%%%%%%%%%%%%%%%%%%%%%%%%%%%%%%%%%%%%%%%%%%%%%%%%%%%%%%%%%%%%%%%%%%%%%%%%%%%%%%%%%%%%%%%%%%%%
\end{thebibliography}
\end{document}